**SARAH PLOSKER and CATHY MATTES**

# Indigenous Beadwork as a Method of Teaching Linear Algebra

**Abstract**

In this work, the authors describe efforts aimed at Indigenizing a second-year linear algebra course at a small liberal arts university in Manitoba, Canada. This is done through an assignment, part hands-on and part written work, that explores the connection between Indigenous beadwork and linear algebra. Our collaboration was perhaps unconventional: Sarah, the first author, is a mathematics professor; while Cathy, the second author, is an associate professor in art history. However, we both had similar goals of putting theory into practice and making positive changes to student learning outcomes in a culturally appropriate way. We situate our work in the context of the current scholarly literature, adding to the important ongoing dialogue on Indigenization of course content and reflecting on the process and outcomes. This transformation of the course curriculum represented an applied approach to immerse Indigenous knowledge and pedagogy into a mathematics classroom. We hope that it may serve as an example of how other educators, particularly in science, technology, engineering, and mathematics (STEM), can integrate Indigenous knowledge-centered pedagogy into their classroom.

## INTRODUCTION

Following the Truth and Reconciliation Commission of Canada Calls to Action (2015), the authors collaborated to create, implement, and assess the impact of an Indigenous beadwork assignment in a second-year undergraduate linear algebra course taught by Sarah Plosker at Brandon University, a small, liberal arts university in southwestern Manitoba, Canada. The idea behind the creation of the course assignment was to see if learning Indigenous beadwork is useful in students' understanding of linear algebra. We felt it was an important first step in aligning our math curriculum with the Truth and Reconciliation Commission of Canada: Calls to Action. The assignment itself was in part a hands-on activity, followed by a reflective written work asking students to connect the mathematical concepts from class with what they learned about beadwork in several in-person beading events.

The main contribution of this work is to provide an example of how mathematics instructors can cultivate cultural relevance in their classroom, putting theory into practice by using beadwork to Indigenize university education. We describe the entire experience in this work, from start to finish. This includes a "pre-project" phase when the authors met and built a trusting relationship. We describe in detail the motivations behind the project, which were three-fold: we hoped the pedagogical approach would be engaging to the students and in turn increase learning, we wanted to respond to the Truth and Reconciliation Commission of Canada Calls to Action, as stated above, and finally this was meant to be a meaningful response to a series of racialized events that took place locally. We situate our work within the larger body of literature, provide a detailed description of our methodology, our experiences, and the specific outcomes of the students' assignment work, which we found overall positive. We conclude with a reflection on the project as a whole, lessons learned, and potential future development if we were to repeat the project in the future.

Local context played a role in how the project was approached. That being said, the ideas behind the project go beyond the local level. In particular, culturally responsive teaching is becoming a worldwide phenomenon—Indigeneity is a concept that goes beyond Canada. We believe a project such as this can show openness, cultural sharing, build trust, and remove some shame for "hesitantly Indigenous" students who may be wrestling with their sense of self. While not meant to be a precise guide to other instructors wishing to do something similar in their classroom, it may serve to spark some ideas of how to implement a similar project in one's own classroom, adapting it to the local context.

We start below by providing some professional details about ourselves, our university, and detailing our motivation for the project.

## CONTEXT

### About the Authors

Sarah Plosker (SP) was at the time an Associate Professor in Mathematics, and held a Canada Research Chair in Quantum Information Theory at Brandon University. She is originally from Regina and started at the University in 2013.

Cathy Mattes (CM) was at the time an Associate Professor of Art History at Brandon University. She taught there from 2002 to 2021. She is now at the University of Winnipeg. She is Southwest Manitoba Métis, and has been beading since she was 20 when she learned from her auntie Jean Baron Ward.

As the University was the center of some major events that led to the Sharing Circle where SP and CM met, and was the institution at which the study took place, we now describe some important details about the University.

**About Brandon University**

Land Acknowledgement (Land Acknowledgement, 2025): Brandon University has campuses on Treaty 1 and Treaty 2 lands, and we are a gathering place for people from many backgrounds and around the world. In this way, we carry on the Indigenous customs of our home in Brandon. We acknowledge Brandon is on shared territory between the Dakota Oyate, the Anishinaabeg, and the National Homeland of the Red River Métis. Today, many other Indigenous people call Brandon their home today, including the Ininew, Anisininewuk, Denesuline, and Inuit.

Please see [About the University] for more information on the University. We provide only a short overview here. Brandon University was founded in 1899 as a Baptist College and became a university in 1967. The Science building, which housed chemistry, physics, botany, zoology, geology, geography, mathematics and computer science, and psychology, was officially opened in 1972. While the Arts and Library building was officially opened in 1961, an Aboriginal Art minor wasn't created until 1993.

In 2019-20, BU has just over 3,000 students, over 300 faculty members, five Faculties (Arts, Education, Graduate Studies, Health Studies, Science), and a School of Music. At present, the Department of Visual and Aboriginal Art is one of the longest running Fine Arts programs in Canada that offers regular Indigenous studio and art history courses. It is one of the only programs in Canada that offers accredited degree majors and minors in Indigenous art, and the only one that requires all of its students to take two Indigenous art history courses, regardless of their Major within Fine Arts.

Our Mission: We promote excellence in teaching, research, creation and scholarship. We educate our students so that they can make a meaningful difference as engaged citizens and leaders. We defend academic freedom and responsibility. We create and disseminate new knowledge. **We embrace cultural diversity and are particularly committed to the education of First Nations, Metis and Inuit people** (emphasis added). We share our expertise and resources with the greater community.

13% of BU students identify as Indigenous (Indigenous Peoples’ Centre, 2025). This number is understood by many to be a low estimate, as students do not always self-identify, for various reasons, on their entrance application to BU. The Indigenous Peoples' Centre on campus offers academic support, transition support (assistance locating e.g. housing, dentist), cultural/spiritual support (cultural events, beading club), and personal and social support (e.g. movie night, soup lunch).

The importance of relationship building with an Indigenous partner such as an Elder, Knowledge Keeper, community member, or in this case, an Indigenous academic professor, in creating a project such as this cannot be overstated; it is the foundation of a reciprocal collaboration and is an important factor in culturally responsive education. In particular, due to historical reasons, Indigenous persons can be skeptical of others calling upon them for their assistance and expertise. There is no formulaic way of forging such relationships, but in general, time, caring, and emphasis on two-way communication are the key components. The collaborative project followed a year of the authors time forming a respectful, trusting, and collegial relationship.

The course was aimed at math and computer science majors, and as such the enrolment is typically small (eight students were enrolled the semester of the study). The beadwork assignment represented an Indigenous pedagogical approach to teaching some of the linear algebra concepts from class. Students were asked to: (1) Attend a beading session in class, hosted by CM, (2) Attend a follow up beading session, open to the community, at the Indigenous Peoples' Center on campus, and (3) Write a short description of how beadwork fits in with topics covered in the math course (roughly 2-3 paragraphs; no more than a page long). Collaborators SP and CM collected data including observations/field notes at in-person beading events both inside and outside the classroom, followed by a qualitative analysis of the written work of the students.

**Motivation for the Project**

We note three main reasons for initiating this project. Firstly, as a response to the Truth and Reconciliation Commission of Canada: Calls to Action (2015). As part of the Indian Residential Schools Settlement Agreement in 2007, the Truth and Reconciliation Commission of Canada was established to facilitate reconciliation among former students, their families, their communities and all Canadians (Truth and Reconciliation Canada Commission of Canada, 2024). The TRC published 94 “Calls to Action”—actionable items meant to redress the legacy of residential schools and advance the process of Canadian reconciliation (Truth and Reconciliation Commission of Canada: Calls to Action, 2015). Our project was an attempt to integrate culturally appropriate topics into the course curriculum. We believe there was a disconnect between our university’s mission statement and student body, and the lack of any Indigenization of course content in Mathematics at the University. We wanted to actively respond to the Calls to Action.

The second main motivation for the collaboration was as a response to a series of white pride events and security risks at the University in December 2017, some of which were reported by the media [Laychuk, 2017]. Both authors attended a Sharing Circle as the situation unfolded, and discussed actions we personally could take going forward. A Sharing Circle is just as it sounds: participants sit in a circle as equals, and a sharing stick is passed around. Only the person with the stick is allowed to talk, and everyone else listens. It is a great way to air grievances, try to find solutions, and connect with each other. Sharing Circles can be highly emotional events, often involving people sharing very personal stories. There may be laughter, crying, and hugging. The fact that SP participated in the Sharing Circle led CM to trust that SP’s intentions were genuine, and this helped pave the way to a trusting, respectful relationship.

Finally, the third reason for initiating the project was the potential for increased learning. As teachers, it becomes obvious quite quickly that not everyone learns best by reading from a textbook. Given our student body, we wanted to create a course assignment based on the Four

R's of First Nations and Higher Education: Respect, Relevance, Reciprocity, and Responsibility (Kirkness & Barnhardt, 1991). In particular, an assignment that would bring relevance to the theoretical aspects we were learning in the course. Mathematical concepts are often quite theoretical, and being able to work things out with one's hands can be useful in solidifying core concepts.

## LITERATURE REVIEW

In recent years, there has been much interest in Indigenizing course content in mathematics, as well as blending arts and mathematics in higher education in Canada. This can be seen both in the academic literature as well as recent relevant national and international conference presentations.

Recent Canadian Mathematical Society Meetings have included sessions on diversity, arts, Indigenization and reconciliation, and their relation to mathematics (The Canadian Mathematical Society, 2019a; The Canadian Mathematical Society, 2019b). These special scientific sessions included a diverse set of talks from mathematics educators across Canada and the USA, with topics ranging from the impact of course design, increasing inclusive and experiential teaching, discussions on active and collaborative learning in mathematics courses, mathematical storytelling, and more. The interest in blending arts and mathematics can also be seen in the USA: Past Joint Mathematics Meetings in the USA have had special sessions on Mathematics and Mathematics Education in Fiber Arts, and now feature an annual Mathematical Art Exhibition (Jensen, 2017; The American Mathematical Society, 2019). While these represent significant contributions to Indigenizing math education content and/or linking math and art, educators would have to be present at the talks; only talk titles and abstracts are available online. In particular, there are no video recordings, transcripts, or conference proceedings. This means that these contributions are of limited use and do not truly serve as references.

An impressive amount of mathematical aspects of Indigenous design are discussed by Doolittle (2017), including geometric transformations with respect to the unified syllabics and birch bark biting, symmetry and groupings in West coast art and totem poles, quill boxes, pueblo pottery, and the geometry, design, and embedded sequences and series in starblankets. This reference was a slideshow given as a presentation by Doolittle at a teachers' conference, which at the time of publication appear to have been taken down, a truly unfortunate aspect of digital work.

Research by Fisher and Mellor (2007) explores group theory, a particular area of mathematics, via the symmetries of a beaded bead. They also explore various patterns and geometric objects, including the sphere, prisms, and pyramids. Later, Fisher and Mellor (2012) focus on flat beading, called angle weaves. This is achieved through the use of mathematical tiling theory (in mathematical terms, "tilings of the plane"). Their works are academic in nature, offering deep mathematical results rather than being experiential math education projects. Belcastro and Yackel (2007) includes chapters discussing quilts, crochet, knitting, cross-stitch, knots, cables, braids, and embroidery, all using mathematical techniques (in particular, geometry and graph theory, among other areas of mathematics). The focus of this book is on the "relationship between mathematics and the fiber arts" and offers methods of introducing the relevant mathematics into the classroom, though it does not specifically involve Indigenous fiber art.

On the topic of Indigenous education, Gaudry and Lorenz (2018) conducted a study based on an anonymous survey of Indigenous academics and their allies. Their work identifies and differentiates three different types of Indigenization happening in the institutional practice of Canadian academia: Indigenous inclusion, reconciliation Indigenization, and decolonial Indigenization. Their work focuses on policy and praxis in the Canadian academy, and is therefore more suited toward administrators versus educators. Ragoonaden and Mueller (2017) conducted a study using student interviews in a general first-year university course at University of British Columbia's Okanagan campus. Their work explored culturally responsive pedagogy and success in first year university courses and programs. Their work showed the importance of circles of learning, peer mentoring, and the student-instructor relationship, and is perhaps the most relevant of the body of literature cited here regarding the impact of culturally responsive pedagogy in university education.

Indigenization efforts in science and mathematics are happening within Canada, the USA, and internationally. Azam and Goodnough (2018) reflected on a self-study of Indigenizing science education in a Science Methods course within the University Education program at a university in Newfoundland (in Atlantic Canada). Azam used a talking circle to learn about and build relationships with the students. Graham (2015) considered Indigenization of the Saskatchewan Mathematics in Canada. Graham discussed various approaches to the goal of linking tipis and mathematics, which evolved into a meta-analysis of their experience relating Indigenous ways of knowing with mathematics. Mirich and Cavey (2015) described mathematics lessons that took the form of moccasin-making to teach the mathematical topics of measurement and area in an elementary school class at a school on a reserve in the USA. Much like the project described herein, the works of Graham and Mirich and Cavey describe specific projects taken on by math educators Indigenizing course content. More generally, Aikenhead (2017) analyses “taken-for-granted notions about school Mathematics” and identifies ways to help enhance school mathematics “in a way that simultaneously promotes both academic achievement and reconciliation”.

On an international level, the mathematics education of Indigenous Mapuche People of Chile is explored by Huencho Ramos (2015), with the idea of implementing culturally relevant teaching activities into mathematics education. The paper took a qualitative approach to this goal, including a five-step process of not only identifying the ethnomathematics knowledge of the Mapuche people, but also analysing how this knowledge is transferred by way of kimches (wise people) in Mapuche communities. A project exploring the revitalization of cultural knowledge through mathematics education, based in New Zealand, was described by Trinick et al. (2015). The paper identified particular challenges, including the lack of local experts, which meant that teachers had to rely on books and other resources, which took much time and effort. Also in New Zealand, three models for culturally responsive teaching in mathematics teacher education are described by Averill et al. (2009). The studies considered the perceptions of preservice and early career teachers in incorporating cultural approaches to mathematics into their teaching. The authors list conditions that they feel are necessary for effective culturally responsive teaching in mathematics: "deep mathematical understanding; effective and open relationships; cultural knowledge; opportunities for flexibility of approach and for implementing change; many accessible and nonthreatening mathematics learning contexts; involvement of a responsive learning community; and, most important, working within a cross-cultural teach." These works all indicate, implicitly or explicitly, that creating a space for culturally relevant mathematics education in one's classroom is no easy task.

Indigenous perspectives on mathematical geometry and space are considered within the context of Papua New Guinea and Australia by Owens (2014). The study found that three principals: language structures, reference lines and points, and measures of space, can be used by teachers to create geometry lesson plans rooted in culture. Abrams et al. (2013) discussed the gap in mathematics and science achievements in public schooling in Australia and the USA, and in STEM careers in the USA. Using a wide variety of sources, they cite the need to create culturally responsive curricula and culturally relevant learning opportunities in K-12 mathematics and science. Our project aims to do precisely this within the university setting.

With SP being a mathematician and CM an Indigenous Art Historian, we felt the need to disseminate our work as there is a notable gap in resources for those outside of critical Indigenous studies. Our project put pedagogies from the different fields of mathematics and Indigenous art into conversation with one another. Although our literature review shows that there are some studies that have engaged in similar questions, there has been a longstanding gap in culturally relevant pedagogical interventions in STEM fields, and we believe our work contributes to this small but growing body of scholarship. Filling the gap in the literature is particularly necessary in light of the Truth and Reconciliation Commission of Canada: Calls to Action, specifically as it relates to Indigenization—integrating Indigenous knowledge, history, and teaching and learning methods into the classroom. Our project is a grounded example that was done in relation to local context and people. Providing culturally relevant learning and Indigenization of course curriculum is of growing global importance, and there is a clear need to document research of this kind. Our aim herein is to contribute to this growing but still vastly underrepresented body of literature, and to exemplify change from the bottom-up.

## METHODOLOGY AND DATA ANALYSIS

An entire class (one hour and 20 minutes) in the second week of the semester was devoted to beadwork. SP wanted the students to learn about beadwork early in the semester, so that she could make reference to various techniques as they were learning mathematical concepts throughout the semester. CM and her Research Assistant came to class and had a large display of beaded items to show the students. CM gave a brief history of beading, descriptions of her own beadwork, and instructions on how to make beaded earrings. Students did hands-on work with the beads, starting their own pair of earrings using “brick stitch”, while CM and her RA gave instructions to the class and helped students one-on-one.

Figure 1 shows examples of the style of earring that the students worked on: a main triangle shape with coloured diamonds in the center of the triangle and strands of beads dangling from the main triangle, using three different colours.

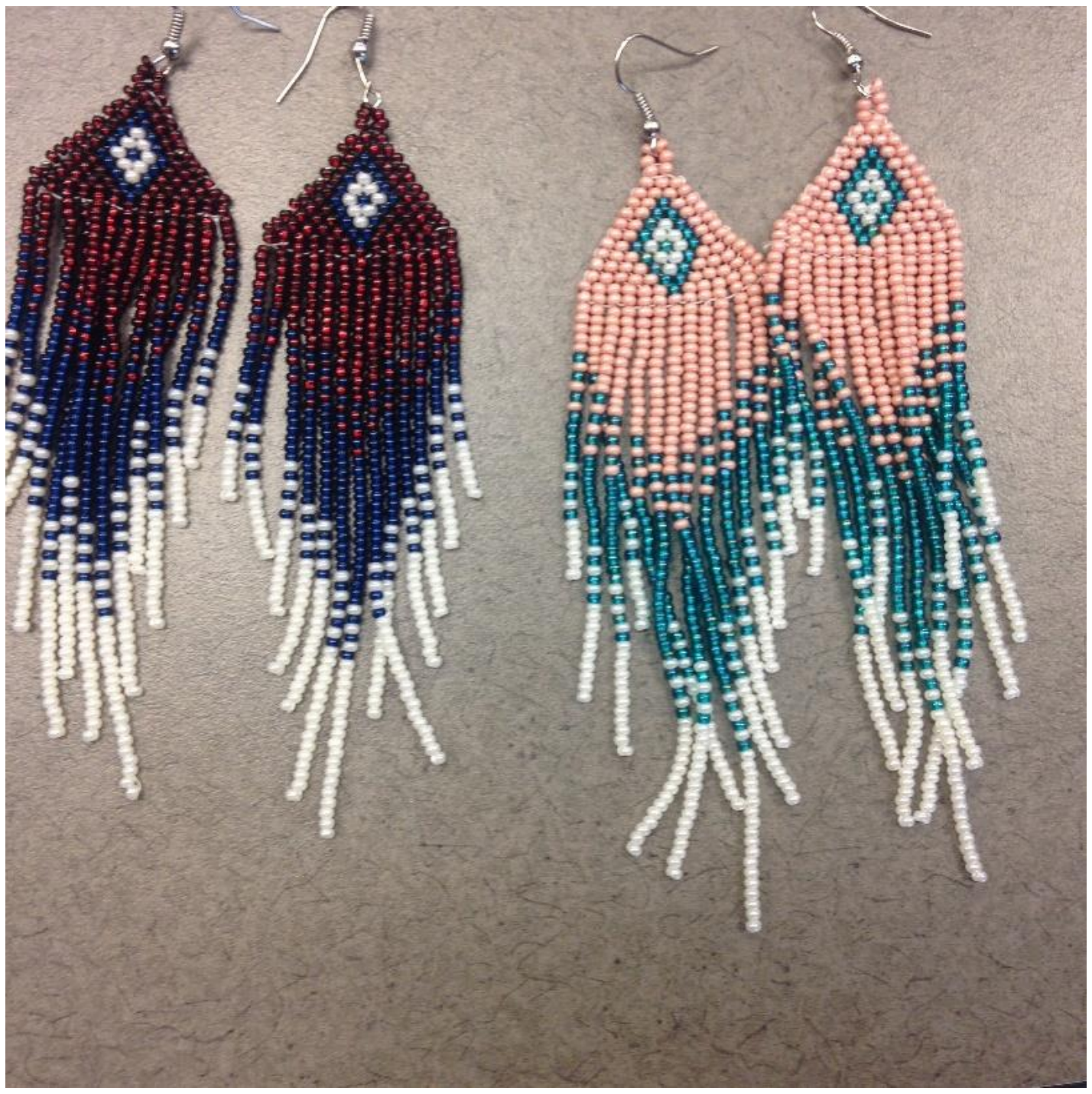

*Figure 1: Examples of earrings of the same style as for the assignment*

Beaded brick stitch earrings are made using a needle and thread. You start with four beads of the same colour on the thread, and then make a "bridge" by stringing the needle through the third and fourth beads. Pulling the thread tight, this changes the line of four beads into a 2-by-2 square of beads. In linear algebraic terms, this is a linear transformation from the four-dimensional real vector space $R^4$ to the set of 2-by-2 matrices $M_2$.

$$\begin{pmatrix} b \\ w \\ g \\ r \end{pmatrix} \rightarrow \begin{pmatrix} w & g \\ b & r \end{pmatrix}$$

*Figure 2: A general 4-dimensional vector being transformed into a 2-by-2 matrix.*

At this point, only one colour of beads has been used, so if we look at the 2-by-2 square we made, it is symmetric in that you can rotate it or flip it over, and it looks the same (this corresponds, for example, to all the letters in Figure 2 being "w", i.e. if all four beads were white). When working with only one colour, some mistakes are not critical. However, if we have different coloured beads, like in Figure 3, then we cannot rotate or flip without changing the pattern.

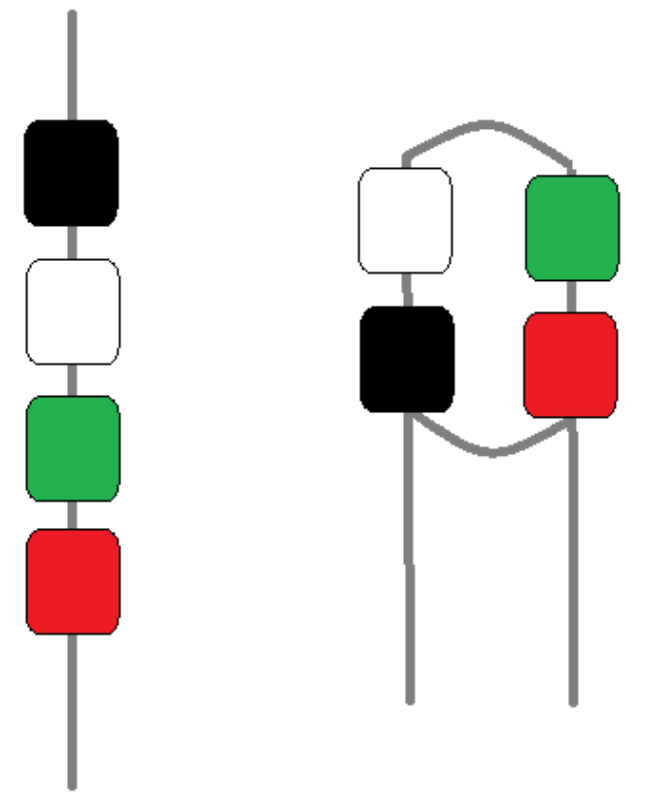

*Figure 3: Four beads (black, white, green, red) forming a column on the left, which then, by threading the needle through the beads in a certain way, becomes a square with two beads in the first row (white, green) and two beads in the second row (black,red). Compare with the mathematical visualization of Figure 2.*

After the initial four beads of the earring, we start adding a second colour of bead, and later, a third colour, so one must be very careful to orient the needle and thread through the beads in the correct direction, otherwise the pattern that emerges from the colours (a diamond shape inside a larger triangle) will not appear correctly. This is a great way to see explicitly what types of matrices are invariant (that is, they don't change) under the different geometric transformations of reflection, rotation, etc. If at any point in making the earrings, you make a mistake, you can carefully undo what you just did, by threading the needle through the beads in the opposite direction. This is, to a mathematician, an invertible function. Many mathematical concepts can be conveyed by this hands-on activity, which brings to life what some students may consider rather stale theoretical content.

Our analysis used qualitative inquiry methods. In the classroom, we relied on observations of the students. Both authors as well as the Research Assistant made mental observations of the students' behaviour (this was not difficult to do, given the small number of students) and discussed our thoughts on several occasions soon after the classroom visit. It was not possible to write field notes during the event as we were all actively involved (with CM and her RA doing the presentation and helping students, and SP taking part in the beading as well as walking around to observe students). Field notes, as such, were written down after the fact; our discussions of our observations and subsequent reflections were recorded in as much accuracy as possible. Participating in this study was voluntary; students submitted a consent form permitting to the use of their written work and observational data in this study.

For the written assignment, both authors read through the submissions separately and then met to discuss the submissions together, reviewing key features of each students' work. Some questions we asked ourselves to guide our analysis were: Did the student clearly articulate the connections that they noticed between the beadwork and course content, as well as their other academic experiences? Did the student relate this connection to their interests in an introspective way? Did the students personalize the assignment writeup in any way?

CM noted several thoughts on Indigenous knowledge sharing in this circumstance: she really appreciated when some of the students started to make diagrams on the chalkboard or drew diagrams on paper. It showed they were engaged, and collaboration is a key aspect of culturally responsive learning. Besides showing their engagement in the material, it also showed active learning strategies which play a key role in Indigenous learning and ways of knowing. Indigenous teachers like CM are used to having those we teach follow our specific directions, as that's how Indigenous knowledge is transmitted. The hardest way is shown first, so that an easier way can be found for each individual after. One will never find an easier way, unless they've learned the way of their knowledge sharer first. Teaching techniques in math are typically one of two ways: start with a simple/low dimensional example, then expand into more difficult examples and a general technique, or start with the general theorem, and then offer a simple example bringing the high-level theorem down to earth.

The second important thought CM noted is regarding the teacher as a learner during the class. This follows the culturally responsive teaching model, being student-centered and eschewing traditional educational practice. This is also quite common with Indigenous knowledge transmission: educators usually place themselves in the position of student, which can be tricky at times. CM appreciated SP's willingness to do so. (During the class when CM and her RA were presenting, SP sat at the desks with the students, facing the front of the classroom). CM finds she is usually challenged when she invites other beadworkers to teach, to learn from them alongside her students, even though she has beaded for most of her adult life. It is challenging as a teacher to not share or interject with one's own thoughts, but important to show respect to a special guest.

The class visit was a condensed learning activity: normally Fine Arts students taking a course on Indigenous art would have multiple lessons prior to making something like earrings. Sorting bags with multiple colours of beads is often a preliminary lesson in patience. In contrast, due to time constraints, we skipped straight to students making earrings for this project. Many students made mistakes and some had to start over multiple times. Some students got out of their desks and started drawing diagrams together on one of the chalkboards in the classroom, trying to get a clearer picture of how the thread was going through the beads to produce the desired outcome. Other students kept to themselves, quietly working away with extreme focus. Some students had to make the mistakes themselves, as they did not believe CM when she explained if they moved the thread in a certain way, they would undo everything. SP notes that this desire of students to check for themselves is somewhat akin to a proof by construction, or reading a mathematical result and trying to prove it oneself rather than reading the proof that the author provides.

Students did not let their frustrations defeat them, and they remained excited throughout the lesson.

The style in which CM and her RA conducted their class visit was very much couched in Indigenous pedagogy, and ways of knowing, being, and doing. CM promoted a deep connection to student learning through an experiential approach, valuing Indigenous knowledge systems in her pedagogy by promoting cultural appreciation and connectedness. The classroom visit was reminiscent of the "trans-systemic" synthesis between Indigenous and Eurocentric knowledge systems discussed in Battiste-Henderson (2021)—the idea of attempting to weave differences and similarities into an overarching method of comprehending.

Beyond our observations of the event, the major data that was collected in this study was the students' written assignments. We discuss these outcomes in the next section.

## STUDENT ASSIGNMENT OUTCOMES

As stated in the previous section, SP wanted the students to learn about beadwork early in the semester, so that she could make reference to various techniques as they were learning mathematical concepts throughout the semester. The idea was that the students would be able to then visualize the beadwork as they learned the corresponding theoretical mathematical concepts throughout the semester. In particular, the first steps of the "brick stitch" transform a line of four beads into a 2-by-2 square of beads (two rows of two beads per row, one stacked below the other). This is an example of what is called a "linear transformation" in mathematics. Students often struggle with the rigorous manner in which linear transformations are taught. The brick stich allowed for a concrete visual example of the abstract notion of a linear transformation. When discussing the mathematical notion of dimension, SP made reference to the fact that the earrings were supposed to lay flat, and pointed out that some mistakes in the beadwork meant that the 2-dimensional objects became 3-dimensional (so, for example, a triangle in the beadwork pattern accidentally became a pyramid). Mistakes were "invertible" (a property of certain linear operations)—if you traced your steps backward carefully, you could undo the mistake you made. These are just some examples of how SP subsequently wove the beadwork demonstration into her pedagogy.

For the written work, SP notes that the "no right answer" written reflection style of the assignment was a stumbling block for some students, since it was so different from their usual math homework. Most of the students approached SP after class or during office hours, in some cases multiple times throughout the semester, to discuss what to write for the assignment. SP managed these discussions similarly to homework help for typical math problems: try to relate the question to topics from class, draw on the students' knowledge from other courses (e.g. the prerequisites), discuss possible approaches, but leave it up to the student to decide in the end how to "solve" the problem. The written work reflected the students' personalities and ways of thinking, with some writeups being more the style of a journal entry, while others were written more like a mathematical textbook—one writeup distilled the student's knowledge of beadwork gained through the beadwork activities into a step-by-step algorithm (many math and computer science students think in terms of algorithms and flowcharts). Each student's work was distinct, but most if not all of the mathematical connections made in class or during discussions outside of class were touched on by at least one student.

We describe and give examples from the written work below. We note that the depth of the written work varied; in some cases it was rather surface level, and in other cases it was quite deep, introspective, and thoughtful. We also note that the written work reflected the students' personalities, and to some degree, whether the student attended a lunch hour or evening beading session at the Indigenous Peoples' Center (students were asked to attend a followup beading activity at the IPC on campus at their convenience during the semester; beading events were held weekly at the IPC at this time). In the evening, the beading events were open to the whole community; grandmothers and artists attended, and there was more interaction, whereas lunch hour sessions tended to be more serious, with students and faculty members working on specific projects. The context and atmosphere of the two types of sessions were quite different.

One student was absent for the beadwork and did not attend the follow-up at the IPC, and as such, that student's work was written rather generically and stood out from the others' written work as lacking introspection as well as connections between the beadwork and their mathematical knowledge. Another student wrote interesting individual sentences but both CM and SP would have liked to have seen additional detail. For example, "the beading session at the Indigenous People's Center…was really more amazing and interesting than the classroom setting". CM and SP both read this in a positive manner, in that it affirmed the importance of the follow-up session at the IPC, but the student did not provide details as to why. Similarly, in the student's mathematical discussion, they state that the beadwork "fits in the matrix representation where you can choose the colors for rows and columns" but did not go into more mathematical detail on matrix representations or provide mathematical examples.

Some students explained the linear transformation from the four-dimensional real vector space $R^4$ to the set of 2-by-2 matrices $M_2$ as illustrated in Figure 2. They went on to explain how certain geometric operations like rotations or flips may or may not change the pattern, depending on exactly what operation is performed, and how many distinct colours are used. Because of symmetry, some mistakes at this stage are not noticeable. However, as was apparent in the classroom demonstration, as we continue with the brick stitch pattern, mistakes become more apparent. They mentioned invertible functions and the connection to carefully undoing a mistake by threading the needle through the beads in the opposite direction.

Many other mathematical concepts can be conveyed by this hands-on activity, which brings to life what some students may consider rather stale theoretical content. Several students went into detail about the dimensions involved in beadwork (i.e., putting the thread through the wrong bead produced a pyramid instead of a triangle). One student discussed how this related to the notions of one-to-one and onto functions between spaces of different dimensions (these functions are ways of transforming one thing to another in some nice mathematical way). One-to-one and onto functions were not discussed explicitly in this manner in class, so this was an analysis beyond simply parroting class discussions.

One student had quite a bit of computer science background, so their written work reflected their discussions with SP regarding planning in advance where each bead will be placed, versus the "painterly style" of beadwork that does not involve pre-planning. The student stated that "free-handing a piece can be rewarding but may not end up how you would like", and explained how a computer program could be helpful "by printing out line by line a series of spaces and stars". This would help artists to avoid gaps of space where a whole bead doesn't quite fit. The student's work indicated that they spent a lot of thought on the intersection of the topics of beading and mathematical programming in computer science, and the authors particularly enjoyed reading their concluding remarks: "At the end of the day a mathematician orders numbers into lines to solve for unknowns in an equation, while an artist places beads onto a string to create something beautiful. Both are simply looking to display rows of information in useful and meaningful ways." The authors note that a matrix displays useful information to a mathematician in the same way that a particular beaded design displays meaning to an artist.

Often mathematicians are rather concrete and like to do things in a logical order. One student, who may have been frustrated by the lack of rigidity in the "rules" for beading earrings, wrote out an algorithm in an attempt to provide more mathematical rigour to the project, almost like a formula for how to make earrings: step one, step two, step three, etc. The student thought of each successive bead that was put on the thread as being labelled 1, 2, 3, … and then would provide instructions such as “Put needle through beads 5+2n and 6+2n through same side as it exits previous sub-step”, where n is the number of desired rows in the “bridge” of the earring. For example, if the bridge was two rows, then this step would indicate to put the needle through beads 9 and 10, in the hole that the needle just exited in the previous step. The student indicated that they wanted to create this algorithm so as to avoid continually asking for instructions. The student also stated some similarities between math and beading, including the idea that if you make a big mistake, then “erase everything (or cut the thread) and start fresh again”. SP enjoyed this student’s writeup as it gives a great deal of insight into a mathematician’s mind, and the potential differences in the ways that a typical math or art student might learn or approach a problem.

One student’s writeup stood out as being incredibly introspective; it was written as though it were a journal or diary entry. CM notes that this parallels Indigenous ways of thinking and storytelling: the student wrote as though they were telling a story. The work contained a great deal of self-reflection and comedy (e.g. MATH as “Mental Abuse Towards Humans”). The student describes their struggles with beading, their frustrations while making mistakes, and the stress of not quite getting it. The student noted the connection with learning math, in particular to their time in grade 12 math in high school, and summarized quite nicely “In dealing with a complex and intricate art form, whether it be beading, or Linear Algebra, it’s very easy to get overwhelmed in just how many different ways there are to mess up. All it takes however, is a few deep breaths and a positive self-affirmation or two to gently nudge yourself back on the right path".

## LESSONS LEARNED

In this collaborative project, we aimed to answer the question of whether or not learning Indigenous beadwork is useful in students' understanding of linear algebra. Incorporating Indigenous ways of learning can help break free from the traditional "chalk talk" that is all too common in academia. Not everyone learns in the same way—typical math courses that are taught straight from the textbook may not be conducive to the majority of the student body. Kinesthetic learning, including tactile, conversation-based, and visual learning (as well as teaching) is more conducive to Indigenous learning and teaching. Indeed, Kovach (2010) writes about the Conversation Method as an Indigenous research methodology that is applicable here, since this was exactly CM's approach to teaching beadwork: conversation-based, without any written instructions. Indigenous scholars such as Simpson (2017) discuss Indigenous knowledge transmission as kinetic and based on relationship-making; a number of the students mentioned positive interactions and building relationships with other beaders at the Indigenous Peoples' Center in the followup in-person event.

As a first attempt at such an exercise, we believe the project was extremely successful, and have a number of ideas for minor revisions to the project in future iterations. In particular, to further blur the divide between art and math, "crits" in studio art courses allow students to obtain feedback on their work and submit their work twice for review by the professor. Allowing students a second attempt, after some initial feedback, would be highly beneficial for the students, as most math and computer science students are not used to this type of assignment.

After the semester ended, Danny Luecke (math and math education faculty at Turtle Mountain Community College, a tribal college in North Dakota) proposed a number of interesting connections to linear algebra topics that SP did not think of during the project (personal communication, 2024). For example, beading on a stretchy material to demonstrate a scaling by pulling and stretching (scaling a vector simply means multiplying each component of the vector by a number; for example, to visualize scaling the vector by two we could stretch a string of beads so that it becomes twice as long as it originally was). This led SP to reflect on other topics, such as adding vectors, which could be done by combining two strings of beads so that the first bead in each of the two strings sit next to each other, the second beads sit next to each other, and so on. The concept of a vector space can be very abstract to a student. One way to visualize a vector space might be to think of the string itself, with its infinite possibilities of combinations of beads that can be put onto it. Besides these missed opportunities, our work did not touch on comparing different beading styles (brick, peyote, lazy, etc) and beading on different surfaces (leggings, moccasins, bonnet, earrings). Such work would likely lead to different interesting mathematical connections, beyond the scope of this article.

**CONCLUSION**

In this paper, we have described a collaborative project whereby the authors created an Indigenous beadwork assignment in an undergraduate mathematics course (specifically, Linear Algebra II). This paper focuses on both the process of the project, including historical and local context, background motivation, and relationship building, and the outcomes of the assignment, including classroom activities, beadwork sessions at the Indigenous Peoples' Center (IPC) on campus, and students' written work. The data collected for this research project included observations by the authors and an RA at the in-person event, as well as a qualitative analysis of the students' written component of the assignment. We also reflected on what went well in the project and things we would change about the project in the future.

The main objectives of this collaboration were to integrate culturally appropriate topics into the course curriculum in an Indigenization and reconciliation effort—to actively respond to the Calls to Action set forth by the Truth and Reconciliation Commission of Canada—and to answer the research question "Does instilling cultural relevance in the classroom curriculum through Indigenous beadwork help with students' understanding of abstract concepts of linear algebra? The impetus for the project came from a series of events that took place locally, leading the authors to meet at a Sharing Circle and ultimately create this project.

Our work puts theory into practice and can be used as a real-world example for other educators wanting to know how to implement culturally responsive teaching “in the field”. A project such as this allows educators to push students out of their comfort zone and gives them permission to think about math in a way that they likely had never thought about before. Our work shows how Indigenous beadwork can contribute to culturally-relevant learning of mathematics, and helps to advance the growing body of work related to culturally-relevant STEM learning discussed in the literature review.

For more information about the project, see https://www.brandonu.ca/research-connection/article/beadwork-and-linear-algebra/ or contact SP at ploskers@brandonu.ca

The project was approved by Brandon University’s Research in Ethics Committee (BUREC) and SP received a certificate of completion of the Tri-Council Policy Statement: Ethical Conduct for Research Involving Humans (TCPS) Course on Research Ethics (CORE) tutorial.